\documentclass{amsart}
\usepackage{graphicx}

\usepackage{amsmath}
\usepackage{amssymb}
\usepackage{amsthm}
\usepackage{calc}
\def\R{{\bf R}}

\def\Diff{\operatorname{Diff}}
\def\ln{\operatorname{ln}}
\def\exp{\operatorname{exp}}
\usepackage[cp1251]{inputenc}
\usepackage[russian]{babel}

\begin{document}

\begin{center}
{
\bf About amenability of  subgroups of the group of
diffeomorphisms of the interval.\\}

 {
 \bf     E.T.Shavgulidze} \\
{\it Department of Mechanics and Mathematics, Moscow State
University, Moscow, 119991 Russia}
 \\
\end{center}

Averaging linear functional on the space continuous functions of
the group of diffeomorphisms of interval is found. Amenability of
several discrete subgroups of the group of diffeomorphisms
$\Diff^3([0,1])$ of interval is prove. In particular, a solution
of the problem of amenability of the Thompson's group $F$ is
given.
\\

{\bf \ \  \ \ \ \ \ \ \ 1. \ \ The main result.}
\\

  Let $\Diff^1_+([0,1])$  be the group of all diffeomorphisms of
class $C^1$ of interval $[0,1]$ that preserve the endpoints of
interval, and  let $\Diff^3_+ ([0,1])$  be the subgroup of
$\Diff^1_+([0,1])$ consisting of all diffeomorphisms  of class
$C^3([0,1])$, and
$$\Diff^3_0 ([0,1])=\{f \in \Diff^3_+ ([0,1]): f'(0)=f'(1)=1 \}.$$
The group $\Diff^1_+([0,1])$ is equipped with the topology
inherited from the space $C^1([0,1])$.

 Let us denote
$ D_n=\{(x_1,...,x_{n-1}):0<x_1<...<x_{n-1}<1\}\subset \R^{n-1}$ \\
and $ x_0=0, \ x_n=1$.

 We say that a subgroup $G$ of
$\Diff^3_0([0,1])$ satisfies condition $(a)$, if

(i) there are a integer $l_n>l_{n-1} \, (l_0=1)$ for any natural
$n$ and a countably additive Borel measure $\eta_n$ on $ D_{l_n}$
such that $\eta_n( D_{l_n})=1$,

(ii) for any positive $\varepsilon$ and any $g \in G$, we can find
natural $N(\varepsilon ,g) $ such that, for any $n>N(\varepsilon
,g) $, it exists a Borel subset $Z_{n,\varepsilon ,g}\subset
D_{l_n}$ that $\eta_n (Z_{n,\varepsilon ,g})>1-\varepsilon$ and
$\max \limits_{1\leq k \leq l_n} (x_k-x_{k-1})<\varepsilon$ for
any $ (x_1,x_2,..., x_{l_n-1})\in Z_{n,\varepsilon ,g}$ where $
x_0=0, \ x_{l_n}=1$,

(iii) $(1-\varepsilon )\eta_n(Y)<\eta_n(gY)< (1+\varepsilon
)\eta_n(Y)$ for any Borel subset $Y\subset Z_{n,\varepsilon ,g}$
where $ gY=\{(g(x_1),g(x_2),..., g(x_{l_n-1})) :(x_1,x_2,...,
x_{l_n-1})\in Y \} $.

For any positive  $\delta <1$, denote by $C_0^{1,\delta}([0,1])$
the set of all functions $f \in C^1([0,1])$ such that $ f(0)=0$
and $  \exists C>0 \ \ \forall t_1,t_2 \in [0,1] \ \
|f'(t_2)-f'(t_1)|<C |t_2-t_1|^\delta.$ Define a Banach structure
on the linear space $C_0^{1,\delta}([0,1])$  by a norm
$$\|f\|_{1,\delta}=|f'(0)|+
\sup \limits_{t_1,t_2 \in
[0,1]}\frac{|f'(t_2)-f'(t_1)|}{|t_2-t_1|^\delta}$$ for any
function $f \in C_0^{1,\delta}([0,1])$.

Let $\Diff^{1,\delta}_+([0,1])=\Diff^1_+([0,1])\bigcap
C_0^{1,\delta}([0,1]).$  It is easy to see that
$\Diff^{1,\delta}_+([0,1])$ is a subgroup of the group
$\Diff^1_+([0,1]).$ The subgroup $\Diff^{1,\delta}_+([0,1])$ is
equipped with the topology inherited from the space
$C_0^{1,\delta}([0,1])$.

Let $C_b(\Diff^{1,\delta}_+([0,1]))$ be the linear space of all
bounded continuous functions on the space
$\Diff^{1,\delta}_+([0,1])$, and let $C_b(\Diff^{1}_+([0,1]))$ be
the linear space of all bounded continuous functions on the space
$\Diff^{1}_+([0,1])$.

Introduce the functions $e_{1,\delta}:
\Diff^{1,\delta}_+([0,1])\to \R, \ \ e_{1,0}:\Diff^{1}_+([0,1])
\to \R$ by setting $e_{1,\delta}(g)=1$ for any $g\in
\Diff^{1,\delta}_+([0,1])$ and $e_{1,0}(f)=1$ \\ for any $f\in
\Diff^{1}_+([0,1])$. Let $F_g(f)=F(g^{-1}\circ f)$ for any $g \in
\Diff^3_0 ([0,1])$,\\ $f\in \Diff^{1,\delta}_+([0,1])$ and $F \in
C_b(\Diff^{1,\delta}_+([0,1]))$.

{\bf Theorem  1.\/} {\it If a  subgroup $G$ of $\Diff^3_0([0,1])$
satisfies condition $(a)$ and a positive $\delta <\frac{1}{2}$
then
there exists a linear functional \\
$L_{\delta}: C_b(\Diff^{1,\delta}_+([0,1]))\to \R$ such that
$L_{\delta}(e_{1,\delta})=1$, $|L_{\delta} (F)|\leq \sup
\limits_{f \in \Diff^{1,\delta}_+([0,1])}|F(f)|$,  $L_{\delta}
(F)\geq 0$ for any nonnegative function $F \in
C_b(\Diff^{1,\delta}_+([0,1]))$, and \\
$L_{\delta} (F_g)=L_{\delta} (F)$ for any $g \in  G$  and $F \in
C_b(\Diff^{1,\delta}_+([0,1]))$. \/}

The restriction of any function of the space $
C_b(\Diff^{1}_+([0,1]))$ on $\Diff^{1,\delta}_+([0,1])$ belongs to
the space $ C_b(\Diff^{1,\delta}_+([0,1]))$. Hence we obtain the
following assertion.

{\bf Corollary  1.1.\/} {\it If a  subgroup $G$ of
$\Diff^3_0([0,1])$ satisfies condition $(a)$ then there exists a
linear functional $L_{0}: C_b(\Diff^{1}_+([0,1]))\to \R$ such that
$L_{0}(e_{1,0})=1$, $|L_{0} (F)|\leq \sup \limits_{f \in
\Diff^{1}_+([0,1])}|F(f)|$, $L_{0} (F)\geq 0$ for any nonnegative
function $F \in C_b(\Diff^{1}_+([0,1]))$, and $L_{0} (F_g)=L_{0}
(F)$ for any $g \in  G$ and $F \in C_b(\Diff^{1}_+([0,1]))$. \/}

 We say that a discrete subgroup $G$ of
$\Diff^3_0([0,1])$ satisfies condition $(b)$, if  there is a such $C>0$ that \\
$\sup \limits_{t \in [0,1]}|\ln(g'_1(t))-\ln(g'_2(t))|\geq C$ for
any $g_1,g_2 \in G, g_1\neq g_2 $.

{\bf Theorem  2.\/} {\it If a discrete subgroup $G$ of
$\Diff^3_0([0,1])$ satisfies conditions $(a),(b)$, then the
subgroup $G$ is amenable. \/}

In [2] \`E.Ghys and V.Sergiescu proved that the Thompson's group
$F$ is isomorphic to a discrete subgroup $G$ of $\Diff^3_0([0,1])$
which satisfies condition $(b)$.

 {\bf Corollary  2.1.\/} {\it The
Thompson's group $F$ is amenable. \/} .
\\

{\bf \ \ \ \ \ \ \ \ \ \ \ \ 2.\ \ \  Proof of Theorem 1.}
\\
\\
 Define the mapping
 $A:\Diff^1_+([0,1])\to C_0([0,1])$ by setting
$$
A(q)(t)=\ln(q'(t))-\ln(q'(0)) \qquad \forall t\in[0,1].
$$

The mapping $A$ is a topological isomorphism between the space
$\Diff^1_+([0,1])$, $C_0([0,1])$ moreover
$$
A^{-1} (\xi) (t)=\frac {\int_0^t e^{\xi (\tau)}d\tau} {\int_0^1
e^{\xi (\tau)}d\tau}.
$$

 Introduce the Wiener measure $w$ on the space $C_0([0,1])$.
Define a Borel  measure $\nu$  on $\Diff^1_+([0,1])$ by setting
$\nu(X)=w(A(X))$ for any Borel subset $X$ of topological space
$\Diff^1_+([0,1])$.

Let $ \delta \in (0, \frac{1}{2})$. It follows from the properties
of Wiener measure $w$ (see [4]) that measure  $\nu$ is
concentrated on the set $E_\delta=\Diff^{1,\delta}_+([0,1])$, i.d.
$\nu(E_\delta)=1$, moreover the Borel subsets of metric space
$E_\delta$ is measurable with respect to the measure $\nu$.

As it was proved in [3], the  measure $\nu$ is quasi-invariant
with respect to the left action of subgroup $\Diff^3_+([0,1])$ on
the group $\Diff^1_+([0,1])$, moreover
\begin{eqnarray}
\nu (g  X)\,=\, \frac 1 {\sqrt{g '(0) g '(1)}} \int_X \,e^{\frac
{g ''(0)} {g '(0)} q'(0)- \frac {g ''(1)} {g '(1)} q'(1) +\int_0^1
S_g (q(t)) (q'(t))^2 dt} \,\nu (d q), \nonumber
\end{eqnarray}
for any  Borel subset $X$  of topological space $\Diff^1_+
([0,1])$,and any
$g \in \Diff^3_+ ([0,1])$,  \\
where $gX=\{ g\circ q : q \in X \}$  and  $ S_g (\tau) = \frac {g
'''(\tau)} {g '(\tau)}- \frac 3 2 (\frac {g ''(\tau)} {g '(\tau)
})^2$ (the Schwartz derivative of function  $g$).

For the proof of Theorem  1 we need the following auxiliary
assertions

{\bf Lemma  1.\/} {\it The following equality is valid
$$\int\limits_{E_\delta}(q'(0))^l\,\nu (d q)=
 \int\limits_{E_\delta}(q'(1))^l\,\nu (d q)$$ for any natural $l$.\/}

Proof. Let $\xi = A(q)$, i.d. $\xi(t)=\ln(q'(t))-\ln(q'(0)) $.
Then
$$
q' (0)=\frac {1} {\int_0^1 e^{\xi (\tau)}d\tau}, \  q' (1)=\frac
{e^{\xi (1)}} {\int_0^1 e^{\xi (\tau)}d\tau}.
$$

Let us take
$$
M_l= \int\limits_{E_\delta}(q'(1))^l\,\nu (d q)
=\int\limits_{\Diff^1_+([0,1])}(q'(1))^l\,\nu (dq)=$$
$$=\int\limits_{C_0([0,1])}(\frac {e^{\xi (1)}} {\int_0^1
e^{\xi (\tau)}d\tau})^l\, w (d \xi)=
 \int\limits_{C_0([0,1])}(\frac {1} {\int_0^1 e^{\xi
(1-\tau)-\xi (1)}d\tau})^l\, w (d \xi)$$

Let $\zeta (t)=\xi (1-t)-\xi (1)$. The Wiener measure $w$ is
invariant with respect to the action $\zeta \longmapsto \xi $,
thus,
$$
M_l= \int\limits_{C_0([0,1])}(\frac {1} {\int_0^1 e^{\zeta
(\tau)}d\tau})^l\,w (d \zeta)=$$
$$=\int\limits_{\Diff^1_+([0,1])}(q'(0))^l\,\nu (dq)=
\int\limits_{E_\delta}(q'(0))^l\,\nu (d q),$$ which implies the
assertion of Lemma 1.

  Introduce the measure $\nu _n=\nu \otimes ...\otimes\nu $ on the space
$E_\delta^n=E_\delta \times ...\times E_\delta $.

Let $c_1=1+M_1+M_2+\int\limits_{E_\delta}(\int_0^1  (q'(t))^2 dt)
\,\nu (d q)$.

For any $r>0$,  $g \in \Diff^3_+([0,1])$,
$\overline{x}=(x_1,...,x_{n-1}) \in D_n$, we write \\
$C_g=1+\max \limits_{0\leq t \leq 1}(|\frac {g''(t)} {g
'(t)}|+(\frac {g''(t)} {g '(t)})^{2}+ |\frac {g'''(t)} {g '(t)}|)$
and
$$X_{r,g,\overline{x}}=\{ (q_1,...,q_n): \ q_1,...,q_n \in E_\delta $$
$$|\sum \limits_{k=1}\limits^{n}[(x_k - x_{k-1})(\frac
{g''(x_{k-1})} {g '(x_{k-1})} q_k'(0)- \frac {g ''(x_k)} {g
'(x_k)} q_k'(1))+$$
$$+(x_k - x_{k-1})^2\int_0^1 S_g (x_{k-1}+(x_k
- x_{k-1})q_k(t)) (q_k'(t))^2 dt ]|\leq 4 c_1 C_g r \}.
$$

{\bf Lemma  2.\/} {\it If $ \epsilon \in (0,1)$, then the
following inequality  is fulfilled \\ $\nu _n(E_\delta^n
\smallsetminus X_{\sqrt[3]{\epsilon},g,\overline{x}})\leq 2
\sqrt[3]{\epsilon}$ for any $g \in  \Diff^3_+([0,1])$,  for any
positive integer  $n$ and
 $\overline{x}=(x_1,...,x_{n-1}) \in D_n$, satisfying
 the inequality \\
 $\max \limits_{1\leq k \leq n}(x_k -
x_{k-1})<\epsilon $.\\
\/}

Proof. Let
$$ f_1 (q_1,...,q_n)=\sum \limits_{k=1}\limits^{n}(x_k - x_{k-1})(\frac
{g''(x_{k-1})} {g '(x_{k-1})} q_k'(0)- \frac {g ''(x_k)} {g
'(x_k)} q_k'(1)).$$

Then
$$
I_1= \int\limits_{E_\delta}...\int\limits_{E_\delta}f_1
(q_1,...,q_n)\,\nu (d q_1)... \nu (d q_n)=M_1 \sum
\limits_{k=1}\limits^{n}(x_k - x_{k-1})(\frac {g''(x_{k-1})} {g
'(x_{k-1})} - \frac {g ''(x_k)} {g '(x_k)} ).$$

As $ |\frac {g''(x_{k-1})} {g '(x_{k-1})} - \frac {g ''(x_k)} {g
'(x_k)}| \leq C_g (x_k - x_{k-1}),$ we have
$$
|I_1| \leq M_1 C_g  \sum \limits_{k=1}\limits^{n}(x_k - x_{k-1})^2
\leq M_1 C_g \epsilon \sum \limits_{k=1}\limits^{n}(x_k -
x_{k-1})=M_1 C_g \epsilon.$$

If $k\neq l$ then
$$\int\limits_{E_\delta}\int\limits_{E_\delta}
(\frac {g''(x_{k-1})} {g '(x_{k-1})}(q_k'(0)-M_1) - \frac {g
''(x_k)} {g '(x_k)} (q_k'(1)-M_1))$$
$$(\frac {g''(x_{l-1})} {g '(x_{l-1})}(q_l'(0)-M_1) - \frac {g
''(x_l)} {g '(x_l)}(q_l'(1)-M_1) )\nu (d q_k) \nu (d q_l)=0,$$
therefore
$$
I_2= \int\limits_{E_\delta}...\int\limits_{E_\delta}(f_1
(q_1,...,q_n)-I_1)^2\,\nu (d q_1)... \nu (d q_n)= $$
$$\sum \limits_{k=1}\limits^{n}(x_k - x_{k-1})^2\int\limits_{E_\delta}[\frac {g''(x_{k-1})} {g
'(x_{k-1})}(q_k'(0)-M_1) - \frac {g ''(x_k)} {g '(x_k)}
(q_k'(1)-M_1)]^2 \nu (d q_k)\leq$$
$$\leq 2 \sum \limits_{k=1}\limits^{n}(x_k - x_{k-1})^2[(\frac {g''(x_{k-1})} {g
'(x_{k-1})})^2 \int\limits_{E_\delta}(q_k'(0)-M_1)^2\nu (d q_k)
+$$
$$+(\frac {g ''(x_k)} {g '(x_k)})^2 \int\limits_{E_\delta}(q_k'(1)-M_1)^2
\nu (d q_k)]=$$
$$= 2 \sum \limits_{k=1}\limits^{n}(x_k - x_{k-1})^2[(\frac {g''(x_{k-1})} {g
'(x_{k-1})})^2 +(\frac {g ''(x_k)} {g '(x_k)})^2]
[\int\limits_{E_\delta}(q_k'(0))^2\nu (d q_k)-(M_1)^2] \leq$$
$$\leq 4 M_2 C_g \sum \limits_{k=1}\limits^{n}(x_k - x_{k-1})^2
 \leq 4 M_2 C_g \epsilon \sum \limits_{k=1}\limits^{n}(x_k
- x_{k-1})=4 M_2 C_g \epsilon.$$

Hence,
$$\nu _n ( \{ (q_1,...,q_n): | f_1
(q_1,...,q_n)-I_1| \geq 2 c_4 C_g \sqrt[3]{\epsilon} \} ) \leq
\frac{I_2}{(2 c_4 C_g \sqrt[3]{\epsilon})^2} \leq \frac{4 M_2 C_g
\epsilon}{(2 c_4 C_g \sqrt[3]{\epsilon})^2} \leq
\sqrt[3]{\epsilon}.
$$
Thus
$$\nu _n ( \{ (q_1,...,q_n): | f_1
(q_1,...,q_n)| \geq 3 c_1 C_g\sqrt[3]{\epsilon} \} ) \leq
\sqrt[3]{\epsilon}.
$$

Let
$$ f_2 (q_1,...,q_n)=\sum \limits_{k=1}\limits^{n}((x_k - x_{k-1})^2\int_0^1 S_g (x_{k-1}+(x_k
- x_{k-1})q_k(t)) (q_k'(t))^2 dt.$$

Then
$$
I_3= \int\limits_{E_\delta}...\int\limits_{E_\delta}|f_2
(q_1,...,q_n)|\,\nu (d q_1)... \nu (d q_n)\leq $$
$$\leq 2 C_g \sum
\limits_{k=1}\limits^{n}(x_k - x_{k-1})^2 \int\limits_{E_\delta}(
\int_0^1 (q_k'(t))^2 dt)\nu (d q_k)\leq$$
$$\leq 2 c_1 C_g \epsilon \sum
\limits_{k=1}\limits^{n}(x_k - x_{k-1})=2 c_1 C_g \epsilon .$$
Thus
$$\nu _n ( \{ (q_1,...,q_n): | f_2
(q_1,...,q_n)| \geq 2 c_1 C_g\sqrt[3]{\epsilon} \} ) \leq
\frac{I_3}{2 c_1 C_g \sqrt[3]{\epsilon}} \leq \frac{2 c_1 C_g
\epsilon}{2 c_1 C_g \sqrt[3]{\epsilon}}= (\sqrt[3]{\epsilon})^2
\leq \sqrt[3]{\epsilon}.
$$

Hence,
$$\nu _n(E_\delta ^n \smallsetminus
X_{\sqrt[3]{\epsilon},g,\overline{x}})= $$
$$=\nu _n ( \{ (q_1,...,q_n): |f_1
(q_1,...,q_n)+ f_2 (q_1,...,q_n)| \geq 4 c_1 C_g\sqrt[3]{\epsilon}
\} )\leq $$
$$\leq \nu _n ( \{ (q_1,...,q_n): |f_1
(q_1,...,q_n)| \geq 2 c_1 C_g\sqrt[3]{\epsilon} \} )+$$
$$+\nu _n ( \{ (q_1,...,q_n): |f_2 (q_1,...,q_n)| \geq 2 c_1 C_g\sqrt[3]{\epsilon}
\} )\leq 2 \sqrt[3]{\epsilon},$$   which implies the assertion of
Lemma 2.

{\bf Lemma 3.\/} {\it For any $g \in  \Diff^3_0([0,1])$,
 $ \epsilon>0$,  there is
$\delta_1 \in (0,1)$ such that the inequality is valid
$$|\prod \limits_{k=1}
\limits^{n}\frac{g(x_k)-g(x_{k-1})}{(x_{k}-x_{k-1})\sqrt{g'(x_k)g'(x_{k-1})}}
-1|\leq \epsilon$$  for any natural $n$ and any
$\overline{x}=(x_1,...,x_{n-1}) \in D_n$ satisfying the inequality
\\ $\max \limits_{1\leq k \leq n}(x_k - x_{k-1})<\delta_1$,  where
$x_0=0$,  $x_n=1$.\/}

Proof. Let $\epsilon \in (0,1)$. Let $C=\max \limits_{t_1, t_2 \in
[0,1]}(1+|\frac{g''(t_{1})}{g'(t_{1})}|+|\frac{g'''(t_{2})}{g'(t_{1})}|)^2$,
$\delta_1 =\frac{1}{400(C+1)}$, $x_{k}'=\frac{x_{k} - x_{k-1}}{2}$
for any $k \ (1\leq k \leq n)$.

There are $x_{k}^*, x_{k}^{**} x_{k}^{***} \in (0,1)$ such that
$$g(x_k)-g(x_{k-1})=g'(x_{k}')(x_k-x_{k-1})+\frac{1}{24}
g'''(x_{k}^*)(x_k-x_{k-1})^3=$$
$$=g'(x_{k}')(x_k-x_{k-1})(1+\frac{g'''(x_{k-1}^*)}{24g'(x_{k}')}
(x_k-x_{k-1})^2),$$
$$g'(x_{k})=g'(x_{k}')+\frac{1}{2}g''(x_{k}')(x_{k}-x_{k-1})+\frac{1}{8}
g'''(x_{k}^{**})(x_{k}-x_{k-1})^2$$
$$=g'(x_{k}')(1+\frac{g''(x_{k}')}{2g'(x_{k}')}
(x_{k-1}-x_{k-2})+\frac{g'''(x_{k}^{**})}{8g'(x_{k}')}
(x_{k-1}-x_{k-2})^2),$$
$$g'(x_{k-1})=g'(x_{k}')-\frac{1}{2}g''(x_{k}')(x_{k}-x_{k-1})+\frac{1}{8}
g'''(x_{k}^{***})(x_{k}-x_{k-1})^2$$
$$=g'(x_{k}')(1-\frac{g''(x_{k}')}{2g'(x_{k}')}
(x_{k-1}-x_{k-2})+\frac{g'''(x_{k}^{***})}{8g'(x_{k}')}
(x_{k-1}-x_{k-2})^2).$$

Hence
$$\frac{g(x_k)-g(x_{k-1})}{(x_{k}-x_{k-1})\sqrt{g'(x_k)g'(x_{k-1})}}=
\frac{1+\lambda _k '(x_k-x_{k-1})} {\sqrt{1+\lambda _k
''(x_k-x_{k-1})}}
$$
where
$$\lambda _k '=\frac{g'''(x_{k-1}^*)}{24g'(x_{k}')}
(x_k-x_{k-1}),$$
$$\lambda _k
''=((\frac{g''(x_{k}')}{2g'(x_{k}')})^2+\frac{g'''(x_{k}^{**})+g'''(x_{k}^{***})}{8g'(x_{k}')})(x_{k}-x_{k-1})+
$$
$$+\frac{g''(x_{k}')(g'''(x_{k}^{***})-g'''(x_{k}^{**}))}{16(g'(x_{k}'))^2}(x_{k}-x_{k-1})^2+
\frac{g'''(x_{k}^{***})g'''(x_{k}^{**})}{64(g'(x_{k}'))^2}(x_{k}-x_{k-1})^3.$$

As $(x_k-x_{k-1})<\delta_1 $, there are $|\lambda _k'|
<C\delta_1<\frac{\epsilon}{100}$, $|\lambda _k''|
<C\delta_1<\frac{\epsilon}{100}$.

We have
$$\sigma=\ln(\prod \limits_{k=1}
\limits^{n}\frac{g(x_k)-g(x_{k-1})}{(x_{k}-x_{k-1})\sqrt{g'(x_k)g'(x_{k-1})}})=
$$
$$=\sum \limits_{k=1}
\limits^{n}(\ln(1+ \lambda _k'(x_k-x_{k-1})-\frac{1}{2}\ln(1+
\lambda _k''(x_k-x_{k-1}))
$$
and
$$|\sigma | \leq 2\sum \limits_{k=1}
\limits^{n}(|\lambda _k'|+|\lambda _k''|)(x_k-x_{k-1})\leq
\frac{\epsilon}{10}\sum \limits_{k=1} \limits^{n}{(x_k-x_{k-1})}=
\frac{\epsilon}{10},
$$
therefore
$$|\prod \limits_{k=1}
\limits^{n}\frac{g(x_k)-g(x_{k-1})}{(x_{k}-x_{k-1})\sqrt{g'(x_k)g'(x_{k-1})}}-1|=
|e^ \sigma   -1| \leq
e^{\frac{\epsilon}{10}}-e^{-\frac{\epsilon}{10}} \leq
\frac{\epsilon}{5}+\frac{\epsilon}{5} < \epsilon,
$$
which implies the assertion of Lemma 3.

Introduce the mapping $Q_n:D_n \times E_\delta^n \to
E_\delta=\Diff^{1,\delta}_+([0,1])$ by setting \\
$f_n\circ (\tilde{l}_n)^{-1}
 =Q_n(x_1,...,x_{n-1},\varphi_1,...,\varphi_{n})$, where
$$f_n(t)=x_{k-1}+(x_k-x_{k-1})\varphi_{k}(n(t-\frac{k-1}{n})),$$
$$\tilde{l}_n (t)=\frac{1}
{x_1-x_0+\sum \limits_{m=2}
\limits^{n}(x_m-x_{m-1})\frac{\varphi_2'(0)\varphi_3'(0)...\varphi_{m}'(0)}
{\varphi_1'(1)\varphi_2'(1)...\varphi_{m-1}'(1)}}\cdot$$
$$\cdot(x_1-x_0+\sum
\limits_{m=2}\limits^{k-1}(x_m-x_{m-1})\frac{\varphi_2'(0)\varphi_3'(0)...\varphi_{m}'(0)}
{\varphi_1'(1)\varphi_2'(1)...\varphi_{m-1}'(1)}+$$
$$+(x_k-x_{k-1})\frac{\varphi_2'(0)\varphi_3'(0)...\varphi_{k}'(0)}
{\varphi_1'(1)\varphi_2'(1)...\varphi_{k-1}'(1)}n(t-\frac{k-1}{n}))
$$
for $t\in [\frac{k-1}{n},\frac{k}{n}]$, $(x_1,...,x_{n-1})\in
D_{n}$, $(\varphi_1,...,\varphi_{n})\in E_\delta^n$.

The function $f=f_n\circ (\tilde{l}_n)^{-1}$ belongs to
$\Diff^{1,\delta}_+([0,1])$, because the left derivation
$$f'( (\tilde{l}_n)^{-1}(\frac{k-1}{n}-0))=n(x_{k-1}-x_{k-2})\varphi_{k-1}(1)\cdot$$
$$\cdot\frac {x_1-x_0+\sum \limits_{m=2}
\limits^{n}(x_m-x_{m-1})\frac{\varphi_2'(0)\varphi_3'(0)...\varphi_{m}'(0)}
{\varphi_1'(1)\varphi_2'(1)...\varphi_{m-1}'(1)}}
{(x_{k-1}-x_{k-2})\frac{\varphi_2'(0)\varphi_3'(0)...\varphi_{k-1}'(0)}
{\varphi_1'(1)\varphi_2'(1)...\varphi_{k-2}'(1)}n }=$$
$$= (x_1-x_0+\sum \limits_{m=2}
\limits^{n}(x_m-x_{m-1})\frac{\varphi_2'(0)\varphi_3'(0)...\varphi_{m}'(0)}
{\varphi_1'(1)\varphi_2'(1)...\varphi_{m-1}'(1)}) \frac
{\varphi_1'(1)\varphi_2'(1)...\varphi_{k-1}'(1)}
{\varphi_2'(0)\varphi_3'(0)...\varphi_{k-1}'(0)}.$$ is equal to
the right derivation
$$f'( (\tilde{l}_n)^{-1}(\frac{k-1}{n}+0))=n(x_k-x_{k-1})\varphi_{k}(0)\cdot$$
$$\cdot\frac {x_1-x_0+\sum \limits_{m=2}
\limits^{n}(x_m-x_{m-1})\frac{\varphi_2'(0)\varphi_3'(0)...\varphi_{m}'(0)}
{\varphi_1'(1)\varphi_2'(1)...\varphi_{m-1}'(1)}}
{(x_k-x_{k-1})\frac{\varphi_2'(0)\varphi_3'(0)...\varphi_{k}'(0)}
{\varphi_1'(1)\varphi_2'(1)...\varphi_{k-1}'(1)}n }=$$
$$= (x_1-x_0+\sum \limits_{m=2}
\limits^{n}(x_m-x_{m-1})\frac{\varphi_2'(0)\varphi_3'(0)...\varphi_{m}'(0)}
{\varphi_1'(1)\varphi_2'(1)...\varphi_{m-1}'(1)}) \frac
{\varphi_1'(1)\varphi_2'(1)...\varphi_{k-1}'(1)}
{\varphi_2'(0)\varphi_3'(0)...\varphi_{k-1}'(0)}.$$

Let a  subgroup $G$ of $\Diff^3_0([0,1])$ satisfies condition
$(a)$.  We write $$L_{\delta , n}(F)=\int\limits_{D_{l_n}}
\int\limits_{E_\delta}...\int\limits_{E_\delta}
F(Q_{l_n}(\bar{x},\varphi_1,...,\varphi_{l_n}))\eta_n (d\bar{x})
\nu (d\varphi_1)...\nu (d\varphi_{l_n})$$ for any function $F \in
C_b(E_\delta)=C_b(\Diff^{1,\delta}_+([0,1]))$.

{\bf Theorem 3.\/} {\it If a  subgroup $G$ of $\Diff^3_0([0,1])$
satisfies condition $(a)$  \\ then $\lim \limits_{n \to
\infty}|L_{\delta , n}(F_g )-L_{\delta , n}(F)|=0$ for any
function $F \in C_b(\Diff^{1,\delta}_+([0,1]))$ and any
diffeomorphism $g \in G$. \/}

Proof. Let $F \in C_b(\Diff^{1,\delta}_+([0,1]))$, $g \in G $, $C=
\sup \limits_{g \in E_\delta}|F(f)|.$

 Let $\epsilon
\in (0,1)$.

It follows from  Lemma 3 that it exists $\delta _1 \in (0,1)$ such
that
$$|\prod \limits_{k=1}
\limits^{n}\frac{g(x_k)-g(x_{k-1})}{(x_{k}-x_{k-1})\sqrt{g'(x_k)g'(x_{k-1})}}
-1|\leq \epsilon$$ for any positive integer  $n$ and for
$\overline{x}=(x_1,...,x_{n-1}) \in D_n$ satisfying the
inequalities $\max \limits_{1\leq k \leq n}(x_k -
x_{k-1})<\delta_1$.

Let us take positive $\epsilon _1$ satisfying the inequalities
$\epsilon _1<\frac{1}{8}\epsilon^3$, $\epsilon _1<\delta_1$, \\
$e^{4c_5 C_g \sqrt[3]{\epsilon _1} }-e^{-4c_4 C_g
\sqrt[3]{\epsilon _1} }<\epsilon$.

 It follows from  Lemma 2 that the inequality is valid $\nu _n(E_\delta ^n
\smallsetminus X_{\sqrt[3]{\epsilon_1},g,\overline{x}})\leq 2
\sqrt[3]{\epsilon_1}\leq \epsilon$ for any positive integer $n$
and for any $\overline{x}=(x_1,...,x_{n-1}) \in D_n$ satisfying
the inequalities   $\max \limits_{1\leq k \leq n}(x_k -
x_{k-1})<\epsilon _1$.

Since the subgroup $G$ satisfies condition $(a)$ we have that

(i) there are a integer $l_n>l_{n-1} \, (l_0=1)$ for any natural
$n$ and a countably additive Borel measure $\eta_n$ on $ D_{l_n}$
such that $\eta_n( D_{l_n})=1$,

(ii)  we can find natural $N(\epsilon _1 ,g) $ such that, for any
$n>N(\epsilon_1 ,g) $, it exists a Borel subset $Z_{n,\epsilon_1
,g}\subset D_{l_n}$ that $\eta_n (Z_{n,\epsilon_1
,g})>1-\epsilon_1$ and $\max \limits_{1\leq k \leq l_n}
(x_k-x_{k-1})<\epsilon_1$ for any $ (x_1,x_2,..., x_{l_n-1})\in
Z_{n,\epsilon_1 ,g}$ where $ x_0=0, \ x_{l_n}=1$,

(iii) $(1-\epsilon_1 )\eta_n(Y)<\eta_n(gY)< (1+\epsilon_1
)\eta_n(Y)$ for any Borel subset $Y\subset Z_{n,\varepsilon ,g}$
where $ gY=\{(g(x_1),g(x_2),..., g(x_{l_n-1})) :(x_1,x_2,...,
x_{l_n-1})\in Y \} $.

Hence it exists the function $\varrho _n :Z_{n,\epsilon_1 ,g} \to
\R$ such that \\ $1-\epsilon_1 \leq \varrho _n (\overline{x})\leq
1+\epsilon_1$ for any $\overline{x} \in Z_{n,\epsilon_1 ,g}$,
$\eta_n(gY)=\int\limits_{Y}\varrho _n (\overline{x}) \eta_n (d
\overline{x})$ for any Borel subset $Y\subset Z_{n,\varepsilon
,g}$.

Let $y_k=g(x_k)$, $\overline{y}=(y_1,...,y_{l_n-1})\in D_{l_n}$,
 \\ $g^{-1}(\overline{y})=(x_1,...,x_{l_n-1})\in
D_{l_n}$. We receive
$$gX_{\sqrt[3]{\epsilon_1},g,g^{-1}(\overline{y})}=\{(\varphi_1,...,\varphi_{l_n})
: (q_1,...,q_{l_n})\in
 X_{\sqrt[3]{\epsilon_1},g,\overline{x}}, \ \ $$
$$ Q_{l_n}(y_1,...,y_{l_n-1},\varphi_1,...,\varphi_{n})=g\circ
(Q_{l_n}(x_1,x_2,...,x_{l_n-1},q_1,...,q_{l_n}))\}
 .$$

It is easy to see that
$\varphi_{k}(t)=\frac{g(x_{k-1}+(x_{k}-x_{k-1})q_k
(t))-g(x_{k-1})}{g(x_{k})-g(x_{k-1})}, $ because
$$\frac{(y_k-y_{k-1})\varphi_2'(0)\varphi_3'(0)...\varphi_{k}'(0)}
{(y_1-y_{0})\varphi_1'(1)\varphi_2'(1)...\varphi_{k-1}'(1)}=
\frac{(x_k-x_{k-1})q_2'(0)q_3'(0)...q_{k}'(0)}
{(x_1-x_{0})q_1'(1)q_2'(1)...q_{k-1}'(1)}$$

We have
$$\int\limits_{gZ_{n,\epsilon _1 ,g}}
 \nu _{l_n}(gX_{\sqrt[3]{\epsilon_1},g,g^{-1}(\overline{y})})
\eta_n (d\overline{y})=
$$
$$=\int\limits_{Z_{n,\epsilon _1 ,g}}
(\int\limits_{X_{\sqrt[3]{\epsilon_1},g,\overline{x}}} \exp(\sum
\limits_{k=1}\limits^{l_n}[(x_k - x_{k-1})(\frac {g''(x_{k-1})} {g
'(x_{k-1})} q_k'(0)- \frac {g ''(x_k)} {g '(x_k)} q_k'(1))+$$
$$+(x_k - x_{k-1})^2\int_0^1 S_g (x_{k-1}+(x_k
- x_{k-1})q_k(t)) (q_k'(t))^2 dt ])\nu (d q_1)... \nu (d q_{l_n}))
$$
$$\varrho _n (\overline{x})\prod \limits_{k=1}
\limits^{l_n}\frac{g(x_k)-g(x_{k-1})}{(x_{k}-x_{k-1})\sqrt{g'(x_k)g'(x_{k-1})}}
\eta_n (d\overline{x})\geq $$
$$\geq(1-\epsilon)^3\int\limits_{Z_{n,\epsilon _1 ,g}}\nu
_{l_n}(X_{\sqrt[3]{\epsilon_1},g,\overline{x}})\eta_n
(d\overline{x})\geq (1-\epsilon)^5.$$

Hence,
$$|L_{\delta ,n} (F_g )-\int\limits_{gZ_{n,\epsilon _1 ,g}}
(
\int\limits_{gX_{\sqrt[3]{\epsilon_1},g,g^{-1}(\overline{y})}}F_g
(Q_{l_n}(\overline{y},\varphi_1,...,\varphi_{l_n}))$$
$$\nu (d\varphi_1)...\nu (d\varphi_{l_n}))
\eta_n (d\overline{y})|\leq C(1-(1-\epsilon)^5)
$$
and
$$|L_{\delta ,n} (F)-\int\limits_{Z_{n,\epsilon _1 ,g}}
( \int\limits_{X_{\sqrt[3]{\epsilon_1},g,\overline{x}}}
F(Q_{l_n}(\overline{x},q_1,...,q_{l_n}))$$
$$\nu (dq_1)...\nu (dq_{l_n}))
\eta_n (d\overline{x})|\leq C(1-(1-\epsilon)^2).
$$

We have
$$|\int\limits_{gZ_{n,\epsilon _1 ,g}}
(
\int\limits_{gX_{\sqrt[3]{\epsilon_1},g,g^{-1}(\overline{y})}}F_g
(Q_{l_n}(\overline{y},\varphi_1,...,\varphi_{l_n}))$$
$$\nu (d\varphi_1)...\nu (d\varphi_{l_n}))
\eta_n (d\overline{y})-
$$
$$-\int\limits_{Z_{n,\epsilon _1 ,g}}(
\int\limits_{X_{\sqrt[3]{\epsilon_1},g,\overline{x}}}
F(Q_{l_n}(\overline{x},q_1,...,q_{l_n}))$$
$$\nu (dq_1)...\nu (dq_{l_n}))
\eta_n (d\overline{x})|\leq
$$
$$\leq\int\limits_{Z_{n,\epsilon _1 ,g}}
(\int\limits_{X_{\sqrt[3]{\epsilon_1},g,\overline{x}}} |\exp(\sum
\limits_{k=1}\limits^{n}[(x_k - x_{k-1})(\frac {g''(x_{k-1})} {g
'(x_{k-1})} q_k'(0)- \frac {g ''(x_k)} {g '(x_k)} q_k'(1))+$$
$$+(x_k - x_{k-1})^2\int_0^1 S_g (x_{k-1}+(x_k
- x_{k-1})q_k(t)) (q_k'(t))^2 dt ])
$$
$$\varrho _n (\overline{x})\prod \limits_{k=1}
\limits^{l_n}\frac{g(x_k)-g(x_{k-1})}{(x_{k}-x_{k-1})\sqrt{g'(x_k)g'(x_{k-1})}}
-1|$$
$$|F(Q_{l_n}(\overline{x},q_1,...,q_{l_n}))|\nu (d q_1)... \nu (d
q_n)) \eta_n (d\overline{x})\leq
$$
$$\leq C\epsilon(2+\epsilon)\int\limits_{Z_{n,\epsilon _1 ,g}}
\nu _n(X_{\sqrt[3]{\epsilon_1},g,\overline{x}})\eta_n
(d\overline{x})\leq C\epsilon(2+\epsilon),$$ which implies the
assertion of Theorem 3.

Define a ultrafilter  $\Im$ on the set positive integers such that
 $\Im$ contains the sets $ \{n,n+1,...\}$ for any positive integer $n$.
We set $L_\delta(F)=\lim \limits_{\Im}L_{\delta , n} ( F)$ for any
function $F \in C_b(E_\delta)$.

Note that the limit always exists because $|L_{\delta , n} (
F)|\leq \sup \limits_{f \in E_\delta}| F(f)|.$

It is easy to see that $L(e_{1,\delta})=1$, $|L_{\delta } (
F)|\leq \sup \limits_{f \in E_\delta}| F(f)|,$ and $L(F)\geq 0$
for any nonnegative function $F \in
C_b(\Diff^{1,\delta}_+([0,1]))$. In turn, Theorem 1 follows from\\
Theorem 3.
\\
\\

{\bf \ \ \ \  \ \ \ \ \ \ \ \ \ \ \ \ 2.\ \ \ Proof of Theorem 2.}
\\
\\
  Let $B(G)$ be the linear space of all bounded  functions on the
group $G$.

Let positive $\delta <\frac{1}{2}$,  let
$$p_\delta(f)=|\ln(f'(0)|+
\sup \limits_{t_1,t_2 \in
[0,1]}\frac{|\ln(f'(t_2))-\ln(f'(t_1))|}{|t_2-t_1|^{\delta}}$$ and
$r(f)=\inf\limits_{h \in G}p_\delta(h^{-1}\circ f)$ for $f \in
\Diff^{1,\delta}_+([0,1])$, $\theta (t)=1-t$ for $t \in [0,1]$ and
$\theta (t)=0$ for $t>1$.

For any fixed $f \in \Diff^{1,\delta}_+([0,1])$,  $C>0$,   the set
of functions \\
$\{\psi : \psi (t)=\ln(g'(t)), \ g \in G, \ \ p_\delta(g\circ
f)<C\} $ contain in a compact subset of the space $C([0,1])$,
therefore it is finite according to condition (а). Hence, we can
define the linear mapping $\pi_{\delta} : B(G) \to
C_b(\Diff^{1,\delta}_+([0,1]))$ by setting
$$\pi_\delta F(f)=\frac{\sum \limits_{h \in G}\theta (p_\delta (h^{-1}\circ f)-r(f))F(h)}
{\sum \limits_{h \in G}\theta (p_\delta (h^{-1}\circ f)-r(f))}.$$

Assign a linear functional $l: B(G)\to \R$ by setting
$l(F)=L_{\delta}(\pi_\delta F)$.

It is easy to see that
$$|l(F)|=|L_{\delta} (\pi_\delta F)|\leq \sup \limits_{f \in
\Diff^{1,\delta}_+([0,1])}|\pi_\delta F(f)|\leq \sup \limits_{g
\in G}|F(g)|,$$ $l(F)\geq 0$ for any nonnegative function $F \in
B(G)$, and $l(e_G)=1$, where $e_G(g)=1$ for all $g \in G$.

Denote by $F_g(h)=F(g^{-1} \circ h)$ for $F \in B(G)$,  $g,h \in
G$.

 We have
$$\pi_\delta F_g(f)=\frac{\sum \limits_{h \in G}\theta (p_\delta(h^{-1}\circ f)-r(f))F(g^{-1}\circ h)}
{\sum \limits_{h \in G}\theta (p_\delta(h^{-1}\circ f)-r(f))}=$$
$$=\frac{\sum \limits_{h \in G}\theta (p_\delta(h^{-1}\circ g\circ
f)-r(g\circ f))F(h)} {\sum \limits_{h \in G}\theta
(p_\delta(h^{-1}\circ g\circ f)-r(g\circ f))}= \pi_\delta F(g\circ
f),$$ hence $l(F_g)=L_{\delta}(\pi_\delta
F_g)=L_{\delta}(\pi_\delta F)=l(F)$,\\ which implies the assertion
of Theorem 2.
\\

 {\bf \ \ \ \ \ \ \ \ \ \ \ \ 3. \ \ \ \ \ \  Proof of Corollary   2.1.}
\\

Let $f_1 (t)=\frac{1}{2}t$ for $0\leq t \leq \frac{1}{2}$,
$f_1(t)=t-\frac{1}{4}$ for $\frac{1}{2}\leq t \leq \frac{3}{4}$,\\
$f_1(t)=2t-1$ for $\frac{3}{4}\leq t \leq 1$ and \\
$f_2(t)=t$ for $0\leq t \leq \frac{1}{2}$,
$f_2(t)=\frac{1}{2}t+\frac{1}{4}$ for $\frac{1}{2}\leq t \leq
\frac{3}{4}$, \\
$f_2(t)=t-\frac{1}{8}$ for $\frac{3}{4}\leq t \leq
\frac{7}{8}$, $f_2(t)=2t-1$ for $\frac{7}{8}\leq t \leq 1$.

The Thompson's group $F$ is generated by $f_1$ and $f_2$.

Denote by  $r_n=1-\frac{1}{2^{n+1}}$ for integer $n\geq 0$ and
$r_{-k}=\frac{1}{2^{k}}$ for integer $k\geq 1$. We have
$f_1(r_n)=r_{n-1}$ for any integer $n$.

The group $F$ act on $ D_n$ by
$f(x_1,...,x_{n-1})=(f(x_1),...,f(x_{n-1}))$ for any $f\in F$,
$(x_1,...,x_{n-1})\in D_n$.

Let$ I_{0}^{k}=\{(r_0,r_1)\}$,
$$
I_{n}^{k}=\{f_2 f_1^{-l_1}f_2 f_1^{l_1-l_2}f_2
...f_1^{l_{n-2}-l_{n-1}}f_2 f_1^{l_{n-1}}(r_0,r_1,...,r_{n+1}): \
0 \leq l_i \leq  \min (k,i)\},
$$
$a_{n,k}=|I_{n}^{k}|$ for $n\geq 0$, $k\geq 1$.

{\bf Lemma 4.\/} {\it $\lim \limits_{ n \to
\infty}\frac{a_{n,k}}{4^{n+1}\cos ^{2n}\frac{\pi}{k+2}} =
(k+2)\sin ^{2}\frac{\pi}{k+2}$.\/}

Proof. It is easy to see that $a_{0,k}=1$, $a_{n,1}=1$,
$a_{n+1,k+1}=\sum \limits_{i=0}\limits^{n}a_{i,k+1}a_{n-i,k}.$

Let $u_{k}(t)=\sum \limits_{n=0}\limits^{\infty}a_{n,k}t^{n}.$

We have $u_{1}(t)=\frac{1}{1-t} $, $u_{k+1}(t)=1+t \
u_{k+1}(t)u_{k}(t)$ or $u_{k+1}(t)=\frac{1}{1-t u_{k}(t)}.$

 Taking $u_{k}(t)=\frac{p_{k-1}(t)}{p_{k}(t)}$, $p_0(t)=1$,
$p_1(t)=1-t$ we find
$$
u_{k+1}(t)=\frac{1}{1-t
\frac{p_{k-1}(t)}{p_{k}(t)}}=\frac{p_{k}(t)}{p_{k}(t)-t
p_{k-1}(t)},
$$
$$p_{k+1}(t)=p_{k}(t)-t
p_{k-1}(t).$$
That means
$$
p_{k}(t)=\frac{1}{2^{k+2}\sqrt{1-4t}}[(1+\sqrt{1-4t})^{k+2}-(1-\sqrt{1-4t})^{k+2}]
$$
or $p_{k}(t)=\prod \limits_{l=1}\limits^{[\frac{k+1}{2}]}(1-4t\cos
^{2}\frac{\pi l}{k+2})$.

Taking $m=[\frac{k+1}{2}]$ we find
$$
u_{k}(t)=\frac{4(k+2)\sin ^{2}\frac{\pi }{k+2}}{1-4t\cos
^{2}\frac{\pi }{k+2} }+...+\frac{4(k+2)\sin ^{2}\frac{\pi
m}{k+2}}{1-4t\cos ^{2}\frac{\pi m}{k+2} }=
$$
$$
=\sum \limits_{n=0}\limits^{\infty}4^{n+1}(k+2)t^{n}(\sin
^{2}\frac{\pi }{k+2} \cos ^{2n}\frac{\pi }{k+2}+...+\sin
^{2}\frac{\pi m}{k+2} \cos ^{2n}\frac{\pi m}{k+2}).
$$

Hence
$$
a_{n,k}=4^{n+1}(k+2)(\sin ^{2}\frac{\pi }{k+2} \cos ^{2n}\frac{\pi
}{k+2}+...+\sin ^{2}\frac{\pi m}{k+2} \cos ^{2n}\frac{\pi m}{k+2})
$$
and $\lim \limits_{ n \to \infty}\frac{a_{n,k}}{4^{n+1}\cos
^{2n}\frac{\pi}{k+2}} = (k+2)\sin ^{2}\frac{\pi}{k+2}$ which
implies the assertion of Lemma 4.

For any integer $l\geq 1$,  $n_1\geq 0$, $n_2\geq 0$,..., $n_l\geq
0$, we write
$$
Y_{l,n_1,n_2,...,n_l}=
\{(r_0,t_{1,1},t_{1,2},...,t_{1,n_1},r_1,t_{2,1},t_{2,2},...,t_{2,n_2}
,r_2,...,
$$
$$r_{l-1},t_{l,1},t_{l,2},...,t_{l,n_l},r_l):\ \
(r_{i-1},t_{i,1},t_{i,2},...,t_{i,n_i},r_i) \in f_1^i
(I_{n_i}^{l-i}), \ \ 0\leq i \leq l \},
$$
$$
Y^{l,n}=\bigcup \limits_{n_1+...+n_l=n,n_1\geq 0,...,n_l\geq 0}
Y_{l,n_1,...,n_l},
$$
$$
Y^{l,n}_0=\bigcup \limits_{n_2+...+n_l=n, n_2\geq 0,...,n_l\geq 0}
Y_{l,0,n_2,...,n_l}.
$$
It is easy to see that $f_2(Y^{l+1,n})= Y^{l,n+1}\setminus
Y^{l,n+1}_0$.

Introduce the mapping $\kappa _n : D_n \to D_{2n}$ by setting
$$
\kappa _n(x_1,x_2,...,
x_{n-1})=(\frac{x_1}{2},x_1,\frac{x_1+x_2}{2},x_2,\frac{x_2+x_3}{2},...,
\frac{x_{n-2}+x_{n-1}}{2},x_{n-1},\frac{x_{n-1}+1}{2}).
$$
Denote by $X^{0,l,n}=\bigcup \limits_{i=0}\limits^{l}\bigcup
\limits_{j=0}\limits^{l-1}f_1^j (Y^{2l-i,n+i})$,
$$
X^{m,l,n}=\kappa _{2^{m-1}(n+2l+2)}(\kappa _{2^{m-2}(n+2l+2)}
(...\kappa _{2(n+2l+2)}(\kappa _{n+2l+2}(X^{0,l,n}))...)).
$$
We have $ X^{m,l,n}=\kappa _{2^{m-1}(n+2l+2)}(X^{m-1,l,n})$ and
$f_1(\kappa _{2^{m-1}(n+2l+2)}(\overline{x}))=\kappa
_{2^{m-1}(n+2l+2)}(f_1(\overline{x}))$, $f_2(\kappa
_{2^{m-1}(n+2l+2)}(\overline{x}))=\kappa
_{2^{m-1}(n+2l+2)}(f_2(\overline{x}))$ for any $\overline{x} \in
X^{m-1,l,n}$.

Also, if $ (t_1,t_2,...,t_{2^{m}(n+2l+2)})$ belongs to $
X^{m,l,n}$ then $$ \{t_1,t_2,...,t_{2^{m}(n+2l+2)}\}\supset \{
\frac{1}{2^m},\frac{2}{2^m},\frac{3}{2^m},...,\frac{2^m-1}{2^m}\}$$
for any $m\geq 1$.

{\bf Lemma 5.\/} {\it For any positive $\varepsilon$, there are
positive integer
$l,n$ such that \\
$\frac{|f_1(X^{m,l,n})\bigcap X^{m,l,n}|}{|X^{m,l,n}|}>1-
\varepsilon$, $\frac{|f_2(X^{m,l,n})\bigcap
X^{m,l,n}|}{|X^{m,l,n}|}>1- \varepsilon$ for any $m\geq 0$ .\/}

Proof. It is sufficient to prove for $m= 0$. Take integer
$l>\frac{4}{\varepsilon}$.

It follows from  Lemma 4 that it exists such integer $n$ that
$\frac{| Y^{2l-i,n+i}_0|}{|Y^{2l-i,n+i}|}<\frac{1}{l^3}$ for any
$0\leq i \leq l$.

As $f_2(Y^{2l-i,n+i})= Y^{l-i-1,n+i+1}\setminus Y^{l-i-1,n+i+1}_0$
we have
$$
(1-\frac{1}{l^2})| Y^{l,n+l}|\leq | Y^{2l,n}|\leq |
Y^{2l-1,n+2}|\leq ...\leq | Y^{l,n+l}|
$$
and $\frac{| Y^{l,n+l}|}{|\bigcup
\limits_{i=0}\limits^{l}Y^{2l-i,n+i}|}<\frac{1}{l}$.

Hence, $\frac{|f_2(X^{0,l,n})\bigcap
X^{0,l,n}|}{|X^{0,l,n}|}>1-\frac{| Y^{l,n+l}|}{|\bigcup
\limits_{i=0}\limits^{l}Y^{2l-i,n+i}|}>1-\frac{1}{l}>1-
\varepsilon$.

As $f_1(\bigcup \limits_{j=0}\limits^{l-1}f_1^j
(Y^{2l-i,n+i}))\bigcap \bigcup \limits_{j=0}\limits^{l-1}f_1^j
(Y^{2l-i,n+i})=\bigcup \limits_{j=1}\limits^{l-1}f_1^j
(Y^{2l-i,n+i})$ \\
we find $\frac{|f_1(X^{0,l,n})\bigcap
X^{0,l,n}|}{|X^{0,l,n}|}=1-\frac{1}{l}>1- \varepsilon$ which
implies the assertion of Lemma 5.

Take a infinite differential function $\psi :\R \to \R$ such that
$\psi (t+1)=\psi (t)+2$, $0< \psi '(t)\leq 3$ for any $t \in \R$,
$\psi '(t)= 3$ for any $t \in [\frac{1}{4},\frac{3}{4}]$, $\psi
(0)= 0$, $\psi (\frac{1}{4})= \frac{1}{4}$, $\psi '(0)= 1$,$\psi
^{(n)}(0)= 0$ for any $n\geq 2$.

For any dyadic rational $r=\frac{k}{2^p} \in (0,1)$, denote $x_r
=\psi ^{-p}(k)$, $x'_r =\psi ^{-p}(k-\frac{1}{4})$, $x_r '' =\psi
^{-p}(k+\frac{1}{4})$, $\phi _r (t)=\psi ^{-p}(k+t)$.

Let $g_1 (t)=\psi ^{-1}(t)$ for $0\leq t \leq
x_{\frac{1}{2}}=\frac{1}{2}$, \\
$g_1(t)=\psi ^{-2}(\psi ^{2}(t)-1)$ for $x_{\frac{1}{2}}\leq t \leq x_{\frac{3}{4}}$,\\
$g_1(t)=\psi (t)-1$ for $x_{\frac{3}{4}}\leq t \leq 1$ and \\
$g_2(t)=t$ for $0 \leq t \leq x_{\frac{1}{2}}$, \\
$g_2(t)=\psi ^{-2}(\psi (t)+1)$ for $x_{\frac{1}{2}}\leq t
\leq x_{\frac{3}{4}}$, \\
$g_2(t)=\psi ^{-3}(\psi ^3(t)-1)$ for $x_{\frac{3}{4}}\leq t \leq
x_{\frac{7}{8}}$, \\
$g_2(t)=\psi (t)-1$ for $x_{\frac{7}{8}}\leq t
\leq 1$.

 In [2] \`E.Ghys and V.Sergiescu proved that the Thompson's
group $F$ is isomorphic to a discrete subgroup $G$ of
$\Diff^3_0([0,1])$ which is generated by $\{g_1, g_2\}$ and
satisfies condition $(b)$.

{\bf Lemma 6.\/} {\it For any dyadic rational $r \in (0,1)$, there
are positive integer $\alpha _1,\alpha _2 ,\beta _1 ,\beta _2$
such that $|\alpha _1|\leq 1$, $|\alpha _2|\leq 1$, $|\beta _1
|\leq 1$, $|\beta _2 |\leq 1$, \\
$g_1(\phi _r (t))=\phi _{f_1 (r)} (\psi ^{\alpha _1}(t))$,
$g_2(\phi _r (t))=\phi _{f_2 (r)} (\psi ^{\beta _1}(t))$,
$g_1(\phi _r (-t))=\phi _{f_1 (r)} (\psi ^{\alpha
_2}(-t))$,\\
$g_2(\phi _r (-t))=\phi _{f_2 (r)} (\psi ^{\beta
_2}(-t))$ for any $t \in [0,\frac{1}{4}]$.\/}

Proof. Let $t \in [0,\frac{1}{4}]$.

If $r =\frac{1}{2}$ we have $f_1 (r)=\frac{1}{4}$, $f_2
(r)=\frac{1}{2}$,
$$g_1(\phi _r (t))=\psi ^{-2}(\psi ^{2}(\psi ^{-1}(t+1))-1)=
\psi ^{-2}(\psi (t)+1)= \phi _{f_1 (r)} (\psi (t)),$$
$$g_2(\phi _r (t))=\psi ^{-2}(\psi ((\psi ^{-1}(t+1)))+1)
=\psi ^{-1}((\psi ^{-1}(t)+1) =\phi _{f_2 (r)} (\psi ^{-1}(t)),$$
$$g_1(\phi _r (-t))=\psi ^{-1}(\psi ^{-1}(-t+1))=\psi ^{-2}(-t+1)=
\phi _{f_1 (r)} (-t),$$
$$g_2(\phi _r (-t))=\psi ^{-1}(-t+1)=\phi _{f_2 (r)} (-t).$$
Hence $\alpha _1 =1,\alpha _2 =0,\beta _1 =-1,\beta _2 =0$.

If $r =\frac{3}{4}$ we have $f_1 (r)=\frac{1}{2}$, $f_2
(r)=\frac{5}{8}$,
$$g_1(\phi _r (t))=\psi (\psi ^{-2}(t+3)-1=
\psi ^{-1}(t+1)= \phi _{f_1 (r)} (t),$$
$$g_2(\phi _r (t))=\psi ^{-3}(\psi ^{3}(\psi ^{-2}(t+3))-1)
=\psi ^{-3}((\psi (t)+5) =\phi _{f_2 (r)} (\psi (t)),$$
$$g_1(\phi _r (-t))=\psi ^{-2}(\psi ^{2}(\psi ^{-2}(-t+3))-1)=\psi ^{-1}(\psi ^{-1}(-t)+1)=
\phi _{f_1 (r)} (\psi ^{-1}(-t)),$$
$$g_2(\phi _r (-t))=\psi ^{-2}(\psi (\psi ^{-2}(-t+3))+1)=
\psi ^{-3}(-t+5)=\phi _{f_2 (r)} (-t).$$
Hence $\alpha _1
=0,\alpha _2 =-1,\beta _1 =1,\beta _2 =0$.

If $r =\frac{7}{8}$ we have  $f_2 (r)=\frac{3}{4}$,
$$g_2(\phi _r (t))=\psi (\psi ^{-3}(t+7))-1
=\psi ^{-2}(t+3) =\phi _{f_2 (r)} (t),$$
$$g_2(\phi _r (-t))=\psi ^{-3}(\psi ^{3} (\psi ^{-3}(-t+7))-1)=
\psi ^{-2}(\psi ^{-1}(-t)+3)=\phi _{f_2 (r)} (-t).$$ Hence
$\beta_1=0,\beta_2=-1$.

If $0<r =\frac{k}{2^p}<\frac{1}{2}$ we have $f_1
(r)=\frac{k}{2^{p+1}}$, $f_2 (r)=\frac{k}{2^p}$,
$$g_1(\phi _r (\pm t))=\psi ^{-1}(\psi ^{-p}(\pm t+k))=\psi ^{-p-1}(\pm t+k)=
\phi _{f_1 (r)} (\pm t),$$
$$g_2(\phi _r (-t))=\psi ^{-p}(\pm t+k)=\phi _{f_2 (r)} (\pm t).$$
Hence $\alpha _1 =\alpha _2 =0,\beta _1 =\beta _2 =0$.

If $\frac{1}{2}<r =\frac{k}{2^p}<\frac{3}{4}$ we have $f_1
(r)=\frac{k-2^{p-2}}{2^p}$, $f_2 (r)=\frac{k+2^{p-1}}{2^{p+1}}$,
$$g_1(\phi _r (\pm t))=\psi ^{-2}(\psi ^{2}(\psi ^{-p}(\pm t+k))-1)=\psi ^{-p}(\pm t+k-2^{p-2})=
\phi _{f_1 (r)} (\pm t),$$
$$g_2(\phi _r (\pm t))=\psi ^{-2}(\psi (\psi ^{-p}(\pm t+k))+1)=
\psi ^{-p-1}(\pm t+k+2^{p-1})=\phi _{f_2 (r)} (\pm t).$$ Hence
$\alpha _1 =\alpha _2 =0,\beta _1 =\beta _2 =0$.

If $\frac{3}{4}<r =\frac{k}{2^p}<1$ we have $f_1
(r)=\frac{k-2^{p-1}}{2^{p-1}}$,
$$g_1(\phi _r (\pm t))=\psi (\psi ^{-p}(\pm t+k))-1=\psi ^{-p+1}(\pm t+k-2^{p-1})=
\phi _{f_1 (r)} (\pm t).$$
 Hence
$\alpha _1 =\alpha _2 =0$.

If $\frac{3}{4}<r =\frac{k}{2^p}<\frac{7}{8}$ we have $f_2
(r)=\frac{k-2^{p-3}}{2^{p}}$,
$$g_2(\phi _r (\pm t))=\psi ^{-3}(\psi ^{3}(\psi ^{-p}(\pm t+k))-1)=\psi ^{-p}(\pm t+k-2^{p-3})=
\phi _{f_2 (r)} (\pm t).$$ Hence $\beta _1 =\beta _2 =0$.

If $\frac{7}{8}<r =\frac{k}{2^p}<1$ we have $f_2
(r)=\frac{k-2^{p-1}}{2^{p-1}}$,
$$g_2(\phi _r (\pm t))=\psi (\psi ^{-p}(\pm t+k))-1=\psi ^{-p+1}(\pm t+k-2^{p-1})=
\phi _{f_2 (r)} (\pm t).$$ Hence $\beta _1 =\beta _2 =0$.

Thus, we prove Lemma 6.

{\bf Lemma 7.\/} {\it For any positive $\varepsilon$,  there are
positive integer $N$ and a finite subset $Z\subset D_N$ such that
$\frac{|g_1(Z)\bigcap Z|}{|Z|}>1- \varepsilon$,
$\frac{|g_2(Z)\bigcap Z|}{|Z|}>1- \varepsilon$,\\
$\max\limits_{1\leq k \leq N} (x_k-x_{k-1})<\varepsilon$ for any
$(x_1,x_2,..., x_{N-1})\in Z$ where $ x_0=0, \ x_{N}=1$ .\/}

Proof. Let $\varepsilon \in (0,1)$.

As $\lim \limits_{ m \to \infty} \sum
\limits_{l=1}\limits^{2^m-1}(x''_{\frac{l}{2^m}}-x'_{\frac{l}{2^m}})=\frac{1}{2}$
it exists such $m\geq 1$ that $\max\limits_{1\leq l \leq 2^m}
(x'_{\frac{l}{2^m}}-x''_{\frac{l-1}{2^m}})<\varepsilon$, where
$x''_0 =\frac{1}{4}$, $x'_1 =\frac{3}{4}$.

By Lemma 5 we find positive integer
$l,n$ such that \\
$\frac{|f_1(X^{m,l,n})\bigcap X^{m,l,n}|}{|X^{m,l,n}|}>1-
\frac{1}{4}\varepsilon$, $\frac{|f_2(X^{m,l,n})\bigcap
X^{m,l,n}|}{|X^{m,l,n}|}>1- \frac{1}{4}\varepsilon$.

Let $k=2^m (n+2l+2)$, $V_{\overline{t}}=\{t_1,t_2,...,t_{k-1}\}$
for any $\overline{t}=(t_1,t_2,...,t_{k-1}) \in X^{m,l,n}$, and
$W=\bigcup \limits_{\overline{t} \in X^{m,l,n}} V_{\overline{t}}$.

Take integer $J>\frac{16(k+1)}{\varepsilon}$. Let
$$
C=\max\limits_{0\leq j \leq J} (\max\limits_{ -\frac{1}{4}\leq x
\leq\frac{1}{4}}(\max\limits_{r \in W} |(\phi _{r} (\psi
^{j}(x)))'|+| (\psi ^{j}(x))'|)).
$$

Take integer $p>\frac{C+1}{\varepsilon}$. Let $N=k(2p+1)$,
$$
Z=\{(\psi ^{j_1}(\frac{1}{4p}),\psi ^{j_1}(\frac{2}{4p}),...,\psi
^{j_1}(\frac{p-1}{4p}),\frac{1}{4},x'_{t_1},\phi _{t_1} (\psi
^{j_2}(-\frac{p-1}{4p})), \phi _{t_1} (\psi
^{j_2}(-\frac{p-2}{4p})),...,
$$
$$
\phi _{t_1} (\psi ^{j_2}(-\frac{1}{4p})),x_{t_1},\phi _{t_1} (\psi
^{j_3}(\frac{1}{4p})),\phi _{t_1} (\psi
^{j_3}(\frac{2}{4p})),...,\phi _{t_1} (\psi
^{j_3}(\frac{p-1}{4p})),x''_{t_1},
$$
$$
x'_{t_2},\phi _{t_2} (\psi ^{j_4}(-\frac{p-1}{4p})), \phi _{t_2}
(\psi ^{j_4}(-\frac{p-2}{4p})),...,\phi _{t_2} (\psi
^{j_4}(-\frac{1}{4p})),
$$
$$
x_{t_2},\phi _{t_2} (\psi ^{j_5}(\frac{1}{4p})),\phi _{t_2} (\psi
^{j_5}(\frac{2}{4p})),...,\phi _{t_2} (\psi
^{j_5}(\frac{p-1}{4p})),x''_{t_2},...,
$$
$$
x'_{t_{k-1}},\phi _{t_{k-1}} (\psi ^{j_{2k-2}}(-\frac{p-1}{4p})),
\phi _{t_{k-1}} (\psi ^{j_{2k-2}}(-\frac{p-2}{4p})),...,\phi
_{t_{k-1}} (\psi ^{j_{2k-2}}(-\frac{1}{4p})),
$$
$$
x_{t_{k-1}},\phi _{t_{k-1}} (\psi ^{j_{2k-1}}(\frac{1}{4p})),\phi
_{t_{k-1}} (\psi ^{j_{2k-1}}(\frac{2}{4p})),...,\phi _{t_{k-1}}
(\psi ^{j_{2k-1}}(\frac{p-1}{4p})),x''_{t_{k-1}},...,
$$
$$
\frac{3}{4},1-\psi ^{j_{2k}}(-\frac{p-1}{4p}), 1-\psi
^{j_{2k}}(-\frac{p-2}{4p}),...,1-\psi ^{j_{2k}}(-\frac{1}{4p})): \
\ 0\leq j_1 \leq J, \ \ 0\leq j_2 \leq J,
$$
$$
\ \ 0\leq j_3 \leq J, \ \ 0\leq j_4 \leq J,\ \ 0\leq j_5 \leq J,
...,\ \ 0\leq j_{2k-2}\leq J,\ \ 0\leq j_{2k-1}\leq J,
$$
$$
0\leq j_{2k}\leq J, \ \ (t_1,t_2,...,t_{k-1}) \in X^{m,l,n} \},
$$
and
$$
Z_i=\{(\psi ^{j_1}(\frac{1}{4p}),\psi
^{j_1}(\frac{2}{4p}),...,\psi
^{j_1}(\frac{p-1}{4p}),\frac{1}{4},x'_{t_1},\phi _{t_1} (\psi
^{j_2}(-\frac{p-1}{4p})), \phi _{t_1} (\psi
^{j_2}(-\frac{p-2}{4p})),...,
$$
$$
\phi _{t_1} (\psi ^{j_2}(-\frac{1}{4p})),x_{t_1},\phi _{t_1} (\psi
^{j_3}(\frac{1}{4p})),\phi _{t_1} (\psi
^{j_3}(\frac{2}{4p})),...,\phi _{t_1} (\psi
^{j_3}(\frac{p-1}{4p})),x''_{t_1},
$$
$$
x'_{t_2},\phi _{t_2} (\psi ^{j_4}(-\frac{p-1}{4p})), \phi _{t_2}
(\psi ^{j_4}(-\frac{p-2}{4p})),...,\phi _{t_2} (\psi
^{j_4}(-\frac{1}{4p})),
$$
$$
x_{t_2},\phi _{t_2} (\psi ^{j_5}(\frac{1}{4p})),\phi _{t_2} (\psi
^{j_5}(\frac{2}{4p})),...,\phi _{t_2} (\psi
^{j_5}(\frac{p-1}{4p})),x''_{t_2},...,
$$
$$
x'_{t_{k-1}},\phi _{t_{k-1}} (\psi ^{j_{2k-2}}(-\frac{p-1}{4p})),
\phi _{t_{k-1}} (\psi ^{j_{2k-2}}(-\frac{p-2}{4p})),...,\phi
_{t_{k-1}} (\psi ^{j_{2k-2}}(-\frac{1}{4p})),
$$
$$
x_{t_{k-1}},\phi _{t_{k-1}} (\psi ^{j_{2k-1}}(\frac{1}{4p})),\phi
_{t_{k-1}} (\psi ^{j_{2k-1}}(\frac{2}{4p})),...,\phi _{t_{k-1}}
(\psi ^{j_{2k-1}}(\frac{p-1}{4p})),x''_{t_{k-1}},...,
$$
$$
\frac{3}{4},1-\psi ^{j_{2k}}(-\frac{p-1}{4p}), 1-\psi
^{j_{2k}}(-\frac{p-2}{4p}),...,1-\psi ^{j_{2k}}(-\frac{1}{4p})): \
\ 1\leq j_1 \leq J-1, \ \ 1\leq j_2 \leq J-1,
$$
$$
1\leq j_3 \leq J-1, \ \ 1\leq j_4 \leq J-1,\ \ 1\leq j_5 \leq J-1,
...,\ \ 1\leq j_{2k-2}\leq J-1,\ \ 1\leq j_{2k-1}\leq J-1,
$$
$$
0\leq j_{2k}\leq J, \ \ (t_1,t_2,...,t_{k-1}) \in
f_i(X^{m,l,n})\bigcap X^{m,l,n} \}
$$
where $i=1,\ \ 2$.

By Lemma 6 we find $g_1 (Z_1)\subset Z$ and $g_2 (Z_2)\subset Z$.
Hence
$$\frac{|g_1 (Z)\bigcap Z|}{|Z|}\leq
\frac{|Z_1|}{|Z|}=\frac{(J-1)^{(2k)}|f_1(X^{m,l,n})\bigcap
X^{m,l,n}|}{(J+1)^{(2k)}| X^{m,l,n}|}>
$$
$$>(1-\frac{4k}{J+1})(1-
\frac{1}{4}\varepsilon )> (1- \frac{1}{4}\varepsilon )^2>1-
\varepsilon ,$$
$$\frac{|g_2 (Z)\bigcap Z|}{|Z|}\leq
\frac{|Z_2|}{|Z|}=\frac{(J-1)^{(2k)}|f_2(X^{m,l,n})\bigcap
X^{m,l,n}|}{(J+1)^{(2k)}| X^{m,l,n}|}> 1- \varepsilon .$$

We have $\phi _{r} (\psi ^{j}(\frac{i}{4p}))-\phi _{r} (\psi
^{j}(\frac{i-1}{4p}))\leq C\frac{1}{4p}<\varepsilon$,    $\psi
^{j}(\frac{i}{4p})-\psi ^{j}(\frac{i-1}{4p})\leq
C\frac{1}{4p}<\varepsilon$ for any $r \in W$, $1\leq i \leq p$,
$1\leq j \leq J$ that means $\max\limits_{1\leq k' \leq N}
(x_{k'}-x_{k'-1})<\varepsilon$ for any $(x_1,x_2,..., x_{N-1})\in
Z$.

Thus, we prove Lemma 7.

In turn, Corollary   2.1 follows from Theorem 2 and Lemma 7.

{\bf Acknowledgements.} The author acknowledge I.K. Babenko,
M.Brin,\\
K.Brown, R. Grigorchuk, V.S. Guba, P. de la Harpe, J.Moore,
M.Sapir, V.Sergiescu, O.G. Smolyanov, A.I. Stern for the
discussions and interest about this work.

 {\bf References}

1. Cannon J.W., Floyd W.J., Parry W.R. "Introductory notes on
Richard\\ Thompson's groups", Enseign Math., vol 42, issue 2
(1996), pages 215--256.

2. Ghys \`E., Sergiescu V. "Sur un groupe remarquable de
diffeomorphismes du cercle", Comment.Math.Helvetici 62 (1987)
185--239.

3.Shavgulidze E.T. "Some Properties of Quasi-Invariant Measures on
Groups of Diffeomorphisms of the Circle",  Russ. J. Math. Physics,
vol 7, issue 4 (2000 ), pages 464--472.

4. Hui-Hsiuhg Kuo "Gaussian measures in Banach space", Lecture
Notes in Mathematics 463, Springer-Verlag 1975.

\end{document}